\documentclass[12pt]{article}       
\usepackage{graphicx, amsbsy, amssymb}

\renewcommand{\baselinestretch}{1.150}

\newcommand{\bu}{\boldsymbol{u}}
\newcommand{\p}{\partial}
\newtheorem{thm}{Theorem}

\newtheorem{pro}{Proposition}

\date{}

\begin{document}

\title{\protect\vspace*{-2cm}\sc\Large New integrable (3+1)-dimensional systems\\ \sc and contact geometry\thanks{This research was performed within the framework of and with financial
support of the OPVK program under project CZ.1.07/2.300/20.0002.
Support from the Ministry of Education,
Youth and Sport of the Czech Republic (M\v{S}MT \v{C}R) under RVO funding for I\v{C}47813059, as well as
from the Grant Agency of the Czech Republic (GA \v{C}R) under grant P201/12/G028,
is also gratefully acknowledged.\newline \hspace*{6mm}This is a post-peer-review, pre-copyedit version of an article published in {\em Letters in Mathematical Physics}. The final
authenticated version is available online at: \href{https://dx.doi.org/10.1007/s11005-017-1013-4}{https://dx.doi.org/10.1007/s11005-017-1013-4}.}
}


\author{{\sc A. Sergyeyev}\\
              Silesian University in Opava, Mathematical Institute,\\ Na Rybn\'\i{}\v{c}ku 1, 74601 Opava, Czech Republic\\
              e-mail: {\tt artur.sergyeyev@math.slu.cz}           
}


\maketitle

\begin{abstract}\protect\vspace*{-1cm}
We introduce a novel systematic construction for integrable (3+1)-dimensional dispersionless systems using nonisospectral Lax pairs that involve contact vector fields. \looseness=-1
%
In particular, we present new large classes of (3+1)-dimensional integrable dispersionless systems associated to the Lax pairs which are polynomial and rational in the spectral parameter.\looseness=-1


\end{abstract}

\section{Introduction}
\label{intro}
Soon after the discovery of the inverse scattering transform
it became clear that integrable nonlinear
systems are of immense significance in modern
mathematics and physics,
cf.\ e.g.\
\cite{ac,al,b-s,bk,cd,dmh,d-book,gph,k,masal,dt-surv}.
Integrable partial differential systems in
four independent variables are, quite naturally, of particular
interest, as, according to general relativity, our spacetime is
four-dimensional. For this reason finding a systematic approach
to the construction of (3+1)-dimensional integrable systems has been
among the most important open problems in modern theory of
integrable systems, see e.g.\  \cite{ac}.\looseness=-1

The overwhelming majority of integrable partial
differential systems in four or more independent variables known to date, cf.\ e.g.\
\cite{ac,bk,d-book,fkk,ms,ms2,tak,zs} and references therein, including the celebrated
(anti-)self-dual Yang--Mills equations and (anti-)self-dual vacuum Einstein equations with vanishing
cosmological constant, can be written as homogeneous first-order
quasilinear, i.e., {\em dispersionless} (also known as {\em hydrodynamic-type}), systems, see e.g.\
\cite{dn,d-book,fkk,z2} and the discussion below for details.

It is therefore natural to look 
for new multidimensional integrable systems which are dispersionless, and we show below that it is indeed possible to construct in a systematic fashion many new (3+1)-dimensional integrable
dispersionless systems using {\em contact Lax pairs}, a novel class of
nonisospectral Lax pairs involving contact vector
fields.

To fix the notation, let $x,y,z,t$ be independent variables, $N$ a natural number (here and below natural numbers mean {\em nonzero nonnegative} integers),
and $u_A$, $A=1,\dots,N$, dependent variables combined into a
vector $\bu=(u_1,\dots,u_N)^{\mathrm{T}}$. From now on the superscript $\mathrm{T}$
indicates the transposed matrix and the subscripts $x,y,z,t,p,u_A$ denote
partial derivatives with respect to the indicated variables. All functions
in the paper
are assumed to be sufficiently smooth for all expressions and computations to make sense;
this can be made rigorous using the language
of differential algebra, see Section 3 in \cite{asp} and references therein for details.\looseness=-1

Recall that a (3+1)-dimensional {\em dispersionless} system is, by definition, a quasilinear homogeneous first-order system, so it can be written in the form
\begin{equation}\label{sys}
A_0 (\bu) \bu_t+A_1(\bu)\bu_x+A_2(\bu)\bu_y+A_3(\bu)\bu_z=0.
\end{equation}
Here 
$A_i=A_i(\bu)$ are $M \times N$ matrices,
$M\geq N$ is a natural number, cf.\ e.g.\
\cite{d-book,fkh,fkk,kg2,ms2}.\looseness=-1

Integrable dispersionless systems in (2+1) dimensions
have the form (\ref{sys}) with $A_3=0$ and $\bu_z=0$.  An
overwhelming majority of integrable systems of this kind,
see e.g.\ \cite{d-book,ms2,z}, can be written as
compatibility conditions of
nonlinear Lax pairs of the form\looseness=-1
\begin{equation}\label{nLax3d}
\psi_y=F(\psi_x,\bu),\quad \psi_t=G(\psi_x,\bu).
\end{equation}

Compatibility of (\ref{nLax3d})
is well known to imply that of a linear
nonisospectral Lax pair written in terms of Hamiltonian
vector fields 
\begin{equation}\label{linlax3d}
\chi_y=\mathcal{X}_f(\chi),
\quad\chi_t =\mathcal{X}_g(\chi).
\end{equation}
Here $\chi=\chi(x,y,t,p)$,
\begin{equation}\label{fg}
f=F(p,\bu),\quad g=G(p,\bu);
\end{equation}
$\mathcal{X}_h$ is a Hamiltonian vector field with the Hamiltonian $h$, i.e.,
\begin{equation}
\mathcal{X}_h=h_p\p_x-h_x\p_p,
\label{hvf}
\end{equation}
and $p$ is an additional independent variable, a {\em variable spectral parameter},
cf.\ e.g.\ \cite{bzm,d-book,ms,z} for details;
recall that $\bu_p\equiv0$ by definition.


The vector field (\ref{hvf}) is two-dimensional, so in order to generalize (\ref{linlax3d}) to (3+1) dimensions we need to replace (\ref{hvf}) by a three-dimensional vector field.\looseness=-1

Contact geometry is a natural odd-dimensional counterpart of the symplectic geometry related to (\ref{hvf}), so we seek to extend the above construction of (2+1)-dimensional integrable dispersionless systems  to (3+1) dimensions by replacing the Hamiltonian vector fields by contact vector fields associated with the contact one-form $dz+pdx$. Such a contact vector field with a contact Hamiltonian $h$ reads
\begin{equation}\label{cvf}
X_h=h_p\p_x-(h_x-p h_z)\p_p+(h-p h_p)\p_z.
\end{equation}
Thus, we propose to replace (\ref{linlax3d}) by
the {\em contact Lax pair} of the form
\begin{equation}\label{clp}
\chi_y=X_f(\chi),\quad \chi_t=X_g(\chi),
\end{equation}
where now $\chi=\chi(x,y,z,t,p)$, $\bu=\bu(x,y,z,t)$, $f=f(p,\bu)$, $g=g(p,\bu)$.
The above system boils down to (\ref{linlax3d}) if $\bu_z=0$ and $\chi_z=0$.

The nonlinear Lax pair (\ref{nLax3d}) also has a counterpart in our setting, a {\em nonlinear contact Lax pair}
\begin{equation}\label{nLax4D}
\psi_y=\psi_z F(\psi_x/\psi_z,\bu),\quad \psi_t=\psi_z G(\psi_x/\psi_z,\bu),
\end{equation}
where $\psi=\psi(x,y,z,t)$ and $F,G$ are related to $f,g$ by the same formula (\ref{fg}) as before.
System (\ref{nLax4D}) boils down to (\ref{nLax3d}) if we put $\bu_z=0$ and $\psi_z=1$.

Unlike the (2+1)-dimensional case, where (\ref{nLax3d}) implies (\ref{linlax3d}) but not vice versa, (\ref{clp}) and (\ref{nLax4D}) have the same compatibility condition
that can be written as a zero-curvature equation (see
Section~\ref{nlp-sec} for the proof):
\begin{thm}\label{zcr-thm}
The following conditions are equivalent:
\begin{center}
\begin{minipage}{0.9\textwidth}
\begin{enumerate}
\item[(i)] a nonisospectral linear contact Lax pair (\ref{clp}) is compatible;
\item[(ii)] nonlinear contact Lax pair (\ref{nLax4D}) with $F,G$ related to $f,g$ from (\ref{clp}) by (\ref{fg}) is compatible;
\looseness=-1
\item[(iii)] the zero-curvature equation
\begin{equation}\label{zcr}
f_t- g_y+\{ f, g\}=0
\end{equation}
holds, where $\{\cdot,\cdot\}$ is the contact bracket:
\begin{equation}\label{cb0}
\begin{array}{rcl}
\{f,g\}
&=&f_{p} g_x-g_{p} f_x-p\left(f_{p} g_z-g_{p} f_z\right)
+ f g_z-g f_z.
\end{array}
\end{equation}
\end{enumerate}
\end{minipage}
\end{center}
\end{thm}


We show below that there exist two fairly broad classes of functions $f(p,\bu)$ and $g(p,\bu)$ such that (\ref{clp}) or, equivalently, (\ref{nLax4D}), gives rise to new non-overdeter\-mined (3+1)-dimensional integrable systems.

For the first of these classes $f$ and $g$ are monic polynomials in $p$ of the form
\begin{equation}\label{polfg}
f=p^{m+1}+
\displaystyle\sum\limits_{i=0}^{m} v_i p^i,\qquad
g=p^{n+1}+\displaystyle\frac{n}{m} v_{m} p^{n}
+\displaystyle\sum\limits_{j=0}^{n-1} w_j p^j,
\end{equation}
where $m$ and $n$
are abitrary natural numbers,
$N=m+n+1$, and the vector $\bu$ of dependent variables reads $\bu=(v_{0},\dots,v_{m},w_{0},\dots,w_{n-1})^{\mathrm{T}}$.

As for the second class, $f$ and $g$ are rational functions of $p$ of the form
\begin{equation}\label{ratfg}
f=\displaystyle\sum\limits_{i=1}^{m}
\frac{a_i}{p-v_i},\quad g=\displaystyle\sum\limits_{j=1}^{n}
\frac{b_j}{p-w_j},
\end{equation}
where $m$ and $n$ are arbitrary natural numbers, so in this case we have $N=2(m+n)$ and
$\bu=(a_1,\dots,a_m,\allowbreak v_1,\dots, v_m,\allowbreak b_1,\dots, b_n,\allowbreak w_1,\dots,w_n)^{\mathrm{T}}$.

The rest of the article is organized as follows. In Section~\ref{mr}
we present new (3+1)-dimensional integrable systems associated with the classes (\ref{polfg}) and (\ref{ratfg}) of Lax functions $f$ and $g$, and the related Lax pairs. In Sections~\ref{lin-lax} and \ref{nlp-sec} we explore the properties of linear and nonlinear contact Lax pairs (\ref{clp}) and (\ref{nLax4D}) as well as the underlying relation to contact geometry, and prove Theorem~\ref{zcr-thm} along the way.
Finally, in Section~\ref{exa} we give
specific examples of new (3+1)-di\-men\-si\-onal integrable
systems related to contact Lax pairs.

\section{New classes of (3+1)-dimensional integrable systems}\label{mr}

As it is common in the literature, we shall say that a
dispersionless
system (\ref{sys}) {\em admits} a linear contact Lax pair (\ref{clp}) (respectively a nonlinear contact Lax pair (\ref{nLax4D})), if (\ref{clp}) (resp.\ (\ref{nLax4D})) is compatible by virtue of (\ref{sys}). Recall that by Theorem~\ref{zcr-thm} admitting (\ref{clp}) is equivalent to admitting (\ref{nLax4D}).

We shall also say that the Lax functions $f$ and $g$ are {\em admissible}
if the compatibility condition for (\ref{clp}) is a system
of the form (\ref{sys}) which is normal in the sense of \cite{o}, i.e., roughly speaking,
it can be transformed into
a system of Cauchy--Kowalevski type by a suitable change of variables.\looseness=-1

There exist plethora of pairs of Lax functions $f$ and $g$ which are admissible
in (2+1) dimensions
in the sense that they give rise to (2+1)-dimensional integrable systems which are normal via (\ref{nLax3d}).
In particular, this includes \cite{z}
monic polynomials and rational functions of $p$ of the form (\ref{ratfg}).
Unfortunately, a complete description of
admissible pairs is missing even in (2+1) dimensions, and hence {\em a fortiori}
in (3+1) dimensions; also, the conditions under which the admissibility in (2+1) dimensions survives in (3+1) dimensions upon passing from (\ref{nLax3d}) to (\ref{nLax4D}) are not known yet.\looseness=-1

Below we explore two broad classes of admissible Lax functions in (3+1) dimensions mentioned in Introduction, (\ref{polfg}) and (\ref{ratfg}).

Recall that for the first of these classes
the functions $f$ and $g$ are monic polynomials in $p$ of the form (\ref{polfg}). The nonlinear contact Lax pair (\ref{nLax4D}) for such $f$ and $g$ 
has the form
\begin{equation}\label{polax}
\begin{array}{l}
\psi_y=\psi_z
\left(\left(\displaystyle\frac{\psi_x}{\psi_z}\right)^{m+1}
+
\displaystyle\sum\limits_{i=0}^{m} v_i \left(\displaystyle\frac{\psi_x}{\psi_z}\right)^i\right),
\\[7mm]
\psi_t=\psi_z\left(\left(\displaystyle\frac{\psi_x}{\psi_z}\right)^{n+1}
\!\!\!\!+\displaystyle\frac{n}{m} v_{m}
\left(\displaystyle\frac{\psi_x}{\psi_z}\right)^{n}
+
\displaystyle\sum\limits_{j=0}^{n-1} w_j
\left(\displaystyle\frac{\psi_x}{\psi_z}\right)^j\right).
\end{array}
\end{equation}

For the sake of brevity put $w_{n}\equiv (n/m)v_m$. 
Then the linear contact Lax pair (\ref{clp}) for our case reads
\begin{equation}\label{clp-pol}
\hspace*{-2mm}
\begin{array}{rcl}
\chi_y&=&\left((m+1) p^{m}+
\displaystyle\sum\limits_{i=1}^{m} i v_i p^{i-1}\right)\chi_x-\left(mp^{m+1}+\displaystyle\sum\limits_{i=0}^{m}(i-1) v_i p^i \right)\chi_z\\[6mm]
&&+\left((v_m)_z p^{m+1}-(v_0)_x+\displaystyle\sum\limits_{i=1}^{m} \left((v_{i-1})_z-(v_i)_x\right) p^i\right)\chi_p,\\[6mm]
\chi_t&=&\left((n+1) p^{n}+
\displaystyle\sum\limits_{j=1}^{n} j w_j p^{j-1}\right)\chi_x-\left(n p^{n+1}+\displaystyle\sum\limits_{j=0}^{n}(j-1) w_j p^j \right)\chi_z\\[6mm]
&&+\Biggl(\displaystyle(w_n)_z p^{n+1}-(w_0)_x+\displaystyle\sum\limits_{j=1}^{n} \left((w_{j-1})_z-(w_j)_x\right) p^j\Biggr)\chi_p.
\end{array}\end{equation}

To simplify writing, assume
that $v_{m+1}=1$, $w_{n+1}=1$,
$v_i=0$ for $i>m+1$ and $i<0$, $w_j=0$ for $j>n+1$ and $j<0$.

Then equating to zero the coefficients at all powers of $p$ in
(\ref{zcr}) for $f$ and $g$ from (\ref{polfg}) yields
%
the system
\begin{equation}\label{gndkp}
\hspace*{-5mm}
\begin{array}{l}
\displaystyle \left(v_k\right)_t-\left(w_k\right)_y
+
\displaystyle\sum\limits_{i=0}^{m+1}
\biggl((k-i-1)w_{k-i} \left(v_i\right)_z-(i-1) v_i \left(w_{k-i}\right)_z\\[7mm]
-(k+1-i)w_{k+1-i} \left(v_i\right)_x
+i v_i \left(w_{k+1-i}\right)_x\biggr)=0,\qquad k=0,\dots,n+m,
\end{array}\hspace{-2mm}
\end{equation}
which is the compatibility condition for both (\ref{polax}) and (\ref{clp-pol}).


The number of equations in (\ref{gndkp}), $n+m+1$, is equal to the number $N$ of dependent variables. 
Moreover, it can be shown that (\ref{gndkp}) can be solved for the
$z$-derivatives $(v_i)_z$ and $(w_j)_z$ for all $i=0,\dots,m$ and $j=0,\dots,n-1$, i.e.,
(\ref{gndkp}) is an evolution system with respect to $z$ in disguise. Thus, the Lax functions $f$ and $g$ from (\ref{polfg}) are admissible for any natural $m$ and $n$.

As $m$ and $n$ in (\ref{polax}) are arbitrary natural numbers,
(\ref{gndkp}) provides an example of a (3+1)-dimensional integrable system with
an arbitrarily large finite number of components.

Now turn to the second of the classes of (3+1)-dimensional integrable systems under study, where $f$ and $g$ are rational functions of the form (\ref{ratfg}),
so the nonlinear contact Lax pair (\ref{nLax4D}) 
reads
\begin{equation}\label{ralaxgen-thm-tr}
\psi_y=\psi_z^2\displaystyle\sum\limits_{i=1}^{m}
\frac{a_i}{\psi_x-v_i\psi_z},
\qquad
\psi_t=\psi_z^2\displaystyle\sum\limits_{j=1}^{n}
\frac{b_j}{\psi_x-w_j\psi_z},
\end{equation}
where $m$ and $n$ are any natural numbers, and
$\bu=(a_1,\dots,a_m,v_1,\dots,v_m,\allowbreak b_1,\dots, b_n,\allowbreak w_1,\dots,w_n)^{\mathrm{T}}$; the convention $w_{n}=(n/m)v_m$ no longer applies.

The linear contact Lax pair (\ref{clp}) in this case is too cumbersome to present its explicit form here.

The compatibility condition for (\ref{ralaxgen-thm-tr}) is
the following (3+1)-dimensional dispersionless
integrable system of $2(m+n)$ equations for $2(m+n)$ unknown functions:
\begin{equation}\label{rls-thm}
\hspace*{-7mm}
\begin{array}{r}
(v_i)_t+\displaystyle\sum\limits_{k=1}^n\left\lbrace
\left(\displaystyle\frac{b_k}{w_k-v_i}\right)_x-\left(\displaystyle\frac{b_k
v_i}{w_k-v_i}\right)_z
+\displaystyle\frac{2 b_k (v_i)_z}{w_k-v_i}\right\rbrace=0,\\[7mm]
(w_j)_y+\displaystyle\sum\limits_{l=1}^m\left\lbrace
-\left(\displaystyle\frac{a_l}{w_j-v_l}\right)_x+\left(\displaystyle\frac{a_l
w_j}{w_j-v_l}\right)_z
-\displaystyle\frac{2 a_l (w_j)_z}{w_j-v_l}\right\rbrace=0,\\[7mm]
(a_i)_t+\displaystyle\sum\limits_{k=1}^n\left\lbrace
\left(\displaystyle\frac{a_i b_k}{(w_k-v_i)^2}\right)_x+\left(\displaystyle\frac{a_i b_k (w_k-2
v_i)}{(w_k-v_i)^2}\right)_z
\right.\\[7mm]\left.
-\displaystyle\frac{3 a_i
(b_k)_z}{w_k-v_i}
+\displaystyle\frac{3 a_i b_k (w_k)_z}{(w_k-v_i)^2}\right\rbrace=0,\\[7mm]
(b_j)_y+\displaystyle\sum\limits_{l=1}^m\left\lbrace
\left(\displaystyle\frac{a_l b_j}{(w_j-v_l)^2}\right)_x+\left(\displaystyle\frac{a_l b_j (w_j-2
v_l)}{(w_j-v_l)^2}\right)_z
\right.\\[7mm]\left.
-\displaystyle\frac{3 a_l
(b_j)_z}{w_j-v_l}
+\displaystyle\frac{3 a_l b_j (w_j)_z}{(w_j-v_l)^2}\right\rbrace=0,
\end{array}
\end{equation}
where $i=1,\dots,m$ and $j=1,\dots,n$.

Passing to new independent variables
$T=t+y$, $Y=t-y$ instead of $t$ and $y$
while keeping $x$, $z$ and all dependent variables intact
turns (\ref{rls-thm})
into a system of Cauchy--Kowalevski type, so
the Lax functions $f$ and $g$ from (\ref{ratfg}) are obviously admissible.\looseness=-1

The number of equations in (\ref{rls-thm})
is equal to the number of
unknown functions, $2(n+m)$. Since $n$ and $m$ are arbitrary natural
numbers by assumption, this system
is, just like (\ref{gndkp}),
an example of a (3+1)-dimensional dispersionless
integrable system with an arbitrarily large finite number of components.\looseness=-1

If $\bu_z=0$ then (\ref{rls-thm}) reduces, modulo the change of signs of
$y$ and $t$, to a (2+1)-dimensional integrable system (5) from
\cite{z}, so (\ref{rls-thm}) can be seen as a (3+1)-dimensional
integrable generalization of the latter system.

In closing note that
the case when $f$ and $g$ 
are rational functions of $p$ in general position, that is,
\begin{equation}\label{ratgp}
f=a_0+\displaystyle\sum\limits_{i=1}^{m}
\frac{a_i}{p-v_i},\quad g=b_0+\displaystyle\sum\limits_{j=1}^{n}
\frac{b_j}{p-w_j},
\end{equation}
can be transformed into (\ref{ralaxgen-thm-tr}) by a suitable change of variables (see \cite{asp} for further details and proof).






\section{Linear contact Lax pairs}\label{lin-lax}

%

To explore the geometric framework related to the linear contact Lax pairs (\ref{clp}),
recall first some basic notions from contact geometry,
see e.g.\ \cite{gm,klr} and references therein.
Namely, consider a smooth 3-manifold $\mathcal{M}$ with local coordinates $x,z,p$ and
contact form $\alpha=dz+pdx$.
A smooth vector field $Y$ on $\mathcal{M}$ is called {\em contact} (with respect to a given contact form $\alpha$)
if there exists a smooth function $\rho$ on $\mathcal{M}$ such that
$L_Y(\alpha)=\rho\alpha$, where $L_Y$ stands for the Lie derivative along $Y$.

A function $h_Y=\alpha(Y)$ is called a {\em contact Hamiltonian}
of $Y$ and we write $Y=X_{h_Y}$ to indicate that a contact vector field $Y$ is uniquely determined by its contact Hamiltonian.
%
For our choice of $\alpha$ and for a given contact Hamiltonian $h$
the associated contact vector field $X_h$ is given by (\ref{cvf}), that is,
\[
X_h
=h_{p} \partial_x+(p
h_z-h_x)\partial_p+(h-p h_{p})\partial_z.
\]

We also have (see e.g.\ \cite{gm,lmp,klr})
the {\em contact bracket} on smooth functions on $\mathcal{M}$
given by (\ref{cb0}), that is,
\[
\{f,g\}
=f_{p} g_x-g_{p} f_x-p\left(f_{p} g_z-g_{p} f_z\right)
+ f g_z-g f_z.
\]
On functions independent of $z$ this bracket boils down to the canonical Poisson bracket in one degree of freedom.
%

It is well known (cf.\ e.g.\ \cite{gm,lmp} and references therein) that the contact bracket is skew-symmetric and obeys the Jacobi identity but not
the Leibniz rule. Instead of the Leibniz rule we have the following generalization thereof:
\[
\{f,g h\}=\{f,g\}h+\{f,h\}g-\{f,1\}gh,
\]
which holds for any smooth functions $f$, $g$, $h$ on $\mathcal{M}$; here $1$ denotes a constant function equal to one.

Thus, the contact bracket (\ref{cb0}) belongs to the class of Jacobi (rather than Poisson) brackets, cf.\ e.g.\ \cite{gm,lmp}, and so this bracket turns the algebra of smooth functions on $\mathcal{M}$ into a Lie
algebra, and, moreover, into a Jacobi algebra, see e.g.\ \cite{gm,lmp} and references therein for
more details on the latter.\looseness=-1


We have (see e.g.\ \cite{gm,lmp,klr})
\begin{equation}\label{hom}
[X_f,X_g]=X_{\{f,g\}},
\end{equation}
for any smooth functions $f$, $g$ on $\mathcal{M}$. Here $[\cdot,\cdot]$ is the standard Lie bracket of vector fields.
Most importantly, the homomorphism $h\mapsto X_h$
from the Lie algebra of smooth functions on $\mathcal{M}$ with respect to the contact bracket to the Lie algebra of smooth vector fields on $\mathcal{M}$
has trivial kernel  (cf.\ e.g.\ \cite{gm,klr}).

The above definitions can be readily applied to the case of the functions
of the form $h=h(p,\bu)$, which depend on $x$ and $z$ via $\bu$. Then for $s\in\{x,y,z,t\}$
we have
\[
h_{s}=\displaystyle\sum_{A=1}^N h_{u_{A}}  (u_A)_{s}.\]

With this in mind, we now turn to the study of linear nonisospectral contact Lax pair (\ref{clp}), that is,
\[
\chi_y=X_f(\chi),\quad \chi_t=X_g(\chi);
\]
recall that $\chi=\chi(x,y,z,t,p)$, $\bu=\bu(x,y,z,t)$, $f=f(p,\bu)$, $g=g(p,\bu)$.

Nonisospectrality means that (\ref{clp})
involves the derivative of $\chi$ with respect to $p$, i.e., $p$ is a variable spectral parameter,
cf.\ e.g.\ \cite{bzm,d-book,ms,z} for details.\looseness=-1

Consider now (cf.\ e.g.\ \cite{ck} and references therein) linear nonisospectral dispersionless Lax pairs which have a more general form than above, namely
\begin{equation}\label{genlinnon}
\chi_y=Q_0\chi_z+Q_1\chi_x+Q_2\chi_p,\quad \chi_t=R_0\chi_z+R_1\chi_z+R_2\chi_p,
\end{equation}
where $\chi=\chi(x,y,z,t,p)$, and the functions $Q_i$ and $R_i$ depend on $\bu$ and on the variable spectral parameter $p$.
It is immediate that (\ref{clp}) is a special case of (\ref{genlinnon}).\looseness=-1

To the best of our knowledge, one can find in the literature only a few integrable dispersionless systems in (3+1) dimensions possessing the Lax pairs of the form (\ref{genlinnon}) which are essentially nonisospectral (i.e., at least one of the functions $Q_2$ and $R_2$ is not identically zero). These systems include {\em inter alia}
the equations for (anti-)self-dual conformal structures \cite{d2}, the Dunajski equation \cite[Section 10.2.3]{d-book}, and certain scalar second-order PDEs from \cite{km}.
However, there appears to be no transformation turning the relevant Lax pairs (\ref{genlinnon})
into the special form (\ref{clp}) for suitable functions $f$ and $g$.
Thus, contact Lax pairs (\ref{clp}) give rise to a genuinely new class of (3+1)-dimensional integrable systems.\looseness=-1

On the other hand, the very fact that contact Lax pairs (\ref{clp}) belong to a broader class of
nonisospectral Lax pairs (\ref{genlinnon}) means that,
at least in principle, 
the contact Lax pairs are amenable to
the inverse scattering transform and to
the dressing method in their various incarnations,
see e.g.\ \cite{dl,kg2,kma,ms,ms2,tt,z2,z} and references therein, and
possibly also to the twistor theory approach, cf.\ e.g.\ \cite{d-book}
and references therein, as well as to some other techniques
like, say, the Hirota method, see \cite{tak} and references therein;
cf.\ also \cite{gh,mma,pre}.\looseness=-1

Let us also point out that
if a system (\ref{sys}) admits a contact Lax pair (\ref{clp})
then we have
\begin{equation}\label{col}
(X_f(\chi))_t-(X_g(\chi))_y=0 
\end{equation}
modulo (\ref{clp}) and its differential consequences.
Substituting a formal expansion of $\chi$ in $p$ into (\ref{col})
yields, upon subsequent splitting with respect to $p$, an infinite hierarchy of
(nonlocal) conservation laws for the associated integrable system
(\ref{sys}), although some of those could be trivial, cf.\ e.g.\
\cite{ac} and Example~1 below.
We anticipate that infinite hierarchies of commuting flows for the
system in question also can be constructed using a formal expansion
of $\chi$ in $p$ in spirit of \cite{tt} and references therein.
\looseness=-1

On a broader note, it would be interesting to explore the structure of symmetries,
conservation laws, symplectic, Hamiltonian and recursion operators, coverings,
and other geometric structures for the systems admitting contact Lax pairs, and to compare
these objects with their counterparts in lower dimensions,
cf.\ e.g.\ \cite{ac,adl,blpla,bsz,dk,kac,dn,d-book,fls,af,bf,h,kv,mawx,my,mokh,o,osw,w}
and references therein.\looseness=-1


In order to gain a better insight into linear contact Lax pairs, consider now the compatibility conditions thereof in more detail.

\begin{pro}\label{zcrcc}Linear contact Lax pair (\ref{clp}) is compatible if and only if the
zero-curvature equation (\ref{zcr}) holds, that is,
\[
f_t- g_y+\{ f, g\}=0.
\]
\end{pro}
\noindent{\em Proof.}
The compatibility condition for (\ref{clp}) obviously reads
\begin{equation}\label{comco}
[\partial_y-X_f,\partial_t-X_g]=0.
\end{equation}

Using (\ref{hom})
we readily see that (\ref{comco}) is equivalent to
\begin{equation}\label{zcrvf}
X_{f_t-g_y+\{f,g\}}=0.
\end{equation}

However, it is immediate from (\ref{cvf}) that
a contact vector field vanishes
if and only if so does its contact Hamiltonian,
hence (\ref{zcrvf}) is equivalent to (\ref{zcr}), and the result follows.

\medskip

In this connection notice that upon
spelling out the derivatives with respect to $x,y,z,t$
in (\ref{zcr}) we find that the latter is equivalent to
the equation
\begin{equation}\label{zcrfull}
\begin{array}{l}
\displaystyle\sum\limits_{A=1}^N \biggl(f_{u_A} (u_A)_t-g_{u_A} (u_A)_{y}+\left(f_p
g_{u_A}-g_p f_{u_A}\right)(u_A)_x
\\[5mm]
+\bigl(\left(f-p f_p\right) g_{u_A}-\left(g- p
g_p\right)f_{u_A}\bigr)(u_A)_z\biggr)=0.
\hspace{-2mm}
\end{array}
\end{equation}
After splitting with respect to $p$ equation (\ref{zcrfull}) gives rise to
a system of the form (\ref{sys}) for $\bu$. Thus,
(\ref{clp}) or, equivalently, (\ref{nLax4D}), whose compatibility condition is (\ref{zcr}),
indeed gives rise 
to integrable systems of
general form (\ref{sys}), i.e., dispersionless (rather than
dispersive) systems, exactly as claimed in Introduction. 

Quite interestingly, existence of a contact Lax pair (\ref{clp}) for a given system (\ref{sys}) entails existence of another linear Lax pair involving the contact bracket:
\begin{pro}\label{clp-pr}
If (\ref{sys}) admits a contact Lax pair (\ref{clp}), it also admits another linear Lax pair:
\begin{equation}\label{cblp}
\begin{array}{l}
\phi_y=\{f,\phi\},\quad
\phi_t=\{g,\phi\},
\end{array}
\end{equation}
where $\phi=\phi(x,y,z,t,p)$.
\end{pro}
\noindent{\em Proof.}
The compatibility condition for (\ref{cblp})
is readily checked to read
\begin{equation}\label{adlaxcc}
\{f_t-g_y+\{f,g\},\phi\}=0,
\end{equation}
so if (\ref{zcr}) holds, then so does (\ref{adlaxcc}), and the result follows.

\medskip

If $\bu_z=0$ and $\phi_z=0$ then (\ref{cblp}) boils down to (\ref{linlax3d}) with $\phi$ replaced by $\chi$. Thus, upon dropping the dependence on $z$ the Lax pairs (\ref{clp}) and (\ref{cblp}) both coalesce into (\ref{linlax3d}), while in (3+1) dimensions they are distinct.

It is easily seen that the solutions of (\ref{cblp})
form a Lie algebra rather than a mere vector space:
\begin{pro}
If $\phi$ and $\tilde\phi$ are two (distinct) solutions of (\ref{cblp}), then their contact bracket $\{\phi,\tilde\phi\}$ also satisfies (\ref{cblp}).
\end{pro}
This property does not seem to have any obvious counterpart for (\ref{clp}).

On the other hand, a straightforward computation readily establishes the following relation among solutions of (\ref{clp}) and (\ref{cblp}):
\begin{pro}
If $\chi$ satisfies (\ref{clp}) and $\phi$ satisfies (\ref{cblp}) then their product $\chi\phi$ again satisfies (\ref{cblp}).

Conversely, if $\phi$ and $\tilde{\phi}$ are two solutions of (\ref{cblp}) and $\phi\not\equiv 0$ then the ratio $\tilde{\phi}/\phi$ satisfies (\ref{clp}).
\end{pro}

In closing note that a solution $\chi$ of (\ref{clp}) has a natural interpretation in the language of dynamical systems.

Namely, suppose that (\ref{clp}) is compatible by virtue of a system (\ref{sys}), and assume for the rest of this section that $\bu$ is a fixed solution of the system (\ref{sys}) under study.

Then $\chi$ is a joint first integral for the following pair of compatible contact dynamical systems on the
contact manifold $\mathcal{M}$,
cf.\ e.g.\ \cite{kt,klr} and references therein,
\begin{equation}
\begin{array}{rcl}
\displaystyle\frac{dx}{dy}&=&-f_p,\\[5mm] \displaystyle\frac{dz}{dy}&=&p f_p-f,\\[5mm] \displaystyle\frac{dp}{dy}&=&-\displaystyle\sum\limits_{A=1}^N f_{u_A} (p (u_A)_z-(u_A)_x),
\end{array}
\label{csf0}
\end{equation}
\begin{equation}
\begin{array}{rcl}
\displaystyle\frac{dx}{dt}&=&-g_p,\\[5mm]
\displaystyle\frac{dz}{dt}&=&p g_p-g,\\[5mm]
\displaystyle\frac{dp}{dt}&=&-\displaystyle\sum\limits_{A=1}^N g_{u_A} (p (u_A)_z-(u_A)_x).
\end{array}\label{csg0}
\end{equation}

Moreover, equations (\ref{nLax4D}) of the associated nonlinear contact Lax pair can be thought of as nonstationary
Hamilton--Jacobi equations for commuting
contact dynamical
systems (\ref{csf0}) and (\ref{csg0}). \looseness=-1
%
%

More precisely, let $H_f=-p_z f(p_x/p_z,\bu)$ and $H_g=-p_z
g(p_x/p_z,\bu)$ be the Hamiltonians with two degrees of freedom
associated with (\ref{nLax4D}). The associated equations of motion
read
\begin{equation}
\begin{array}{rclcrcl}
\displaystyle\frac{dx}{dy}&=&(H_f)_{p_x}=-f_{\tilde{p}},&\quad& 
\displaystyle\frac{dz}{dy}&=&(H_f)_{p_z}=\tilde{p} f_{\tilde p}-f,\\[5mm]
\displaystyle\frac{dp_x}{dy}&=&-(H_f)_x=p_z f_x,&\quad& 
\displaystyle\frac{dp_z}{dy}&=&-(H_f)_{z}=p_z f_z,\end{array}
\label{pxpzsysf0}
\end{equation}
\begin{equation}
\begin{array}{rclcrcl}
\displaystyle\frac{dx}{dt}&=&(H_g)_{p_x}=-g_{\tilde{p}}, &\quad& 
\displaystyle\frac{dz}{dt}&=&(H_g)_{p_z}=\tilde{p} g_{\tilde p}-g,\\[5mm]
\displaystyle\frac{dp_x}{dt}&=&-(H_g)_x=p_z g_x, &\quad& 
\displaystyle\frac{dp_z}{dt}&=&-(H_g)_{z}=p_z g_z,
\end{array}\label{pxpzsysg0}
\end{equation}
where $\tilde{p}=p_x/p_z$.

Equations (\ref{nLax4D}) are nothing but the nonstationary Hamilton--Jacobi
equations for these dynamical systems. Moreover, the flows  of
(\ref{pxpzsysf0}) and (\ref{pxpzsysg0}) commute if the system (\ref{nLax4D})
is compatible.\looseness=-1

It readily follows that we have
\[
\displaystyle\frac{d\tilde{p}}{dy}=f_x-\tilde{p} f_z,\quad \displaystyle\frac{d\tilde{p}}{dt}=g_x-\tilde{p}
g_z,
\]
so we can rewrite (\ref{pxpzsysf0}) and (\ref{pxpzsysg0}) as
\begin{equation}\label{pxpzsysf}
\begin{array}{rclcrcl}
\displaystyle\frac{dx}{dy}&=&-f_{\tilde{p}}, &\quad& \displaystyle\frac{dz}{dy}&=&\tilde{p} f_{\tilde p}-f, \\[5mm]
\displaystyle\frac{d\tilde{p}}{dy}&=&f_x-\tilde{p} f_z, &\quad& \displaystyle\frac{dp_z}{dy}&=&p_z f_z,
\end{array}
\end{equation}
\begin{equation}\label{pxpzsysg}
\begin{array}{rclcrcl}
\displaystyle\frac{dx}{dt}&=&-g_{\tilde{p}}, &\quad& \displaystyle\frac{dz}{dt}&=&\tilde{p} g_{\tilde p}-g, \\[5mm]
\displaystyle\frac{d\tilde{p}}{dt}&=&g_x-\tilde{p} g_z, &\quad& \displaystyle\frac{dp_z}{dt}&=&p_z
g_z.
\end{array}
\end{equation}

It is now clear that if we identify $\tilde{p}$ with $p$ then
(\ref{pxpzsysf}) and (\ref{pxpzsysg}) are nothing but (\ref{csf0}) and (\ref{csg0})
supplemented by the equations for $p_z$,
\[
\begin{array}{l}
\displaystyle\frac{dp_z}{dy}=p_z f_z,\quad 
\displaystyle\frac{dp_z}{dt}=p_z g_z, 
\end{array}
\]
which are straightforward to solve once (\ref{pxpzsysf}) and
(\ref{pxpzsysg}) are solved. For this very reason it is natural to
consider (\ref{nLax4D}) as the Hamilton--Jacobi equations not just for
(\ref{pxpzsysf0}) and (\ref{pxpzsysg0}) but also for (\ref{csf0}) and
(\ref{csg0}).

\section{Nonlinear contact Lax pairs}\label{nlp-sec}

We start with noticing that the compatibility condition for a nonlinear contact Lax pair (\ref{nLax4D}) is the same as for its linear counterpart (\ref{clp}):
\begin{pro}\label{zcr-nlp}Nonlinear contact Lax pair (\ref{nLax4D}) is compatible if and only if the
zero-curvature equation (\ref{zcr}) holds, that is,
\[
f_t- g_y+\{ f, g\}=0,
\]
where $f=F(p,\bu)$ and $g=G(p,\bu)$.
\end{pro}
\noindent{\em Proof.} Spelling out the compatibility condition $(\psi_y)_t-(\psi_t)_y=0$ for (\ref{nLax4D})
and taking into account (\ref{nLax4D}) itself
reveals that the said condition is equivalent to the equation
\begin{equation}\label{zcrnlp}
\left(f_t- g_y+\{f, g\}\right)_{p=\psi_x/\psi_z}=0,
\end{equation}
where $f=F(p,\bu)$ and $g=G(p,\bu)$.
Equation (\ref{zcrnlp}) is nothing but (\ref{zcr}) modulo the substitution $p=\psi_x/\psi_z$,
and the result follows.

\medskip

Theorem~\ref{zcr-thm} now readily follows from Propositions~\ref{zcrcc} and \ref{zcr-nlp}.
%
Indeed, by  Proposition~\ref{zcrcc} 
(i) implies (iii) while by Proposition~\ref{zcr-nlp} (ii) implies (iii). On the other hand, it is immediate from Propositions~\ref{zcrcc} and \ref{zcr-nlp} that (iii) implies both (i) and (ii), and Theorem~\ref{zcr-thm} is proved.






Even though Theorem~\ref{zcr-thm} establishes equivalence of linear and nonlinear contact Lax pairs, 
the latter are of interest on their own right.

First of all, some of the nonlinear contact Lax pairs could admit dispersive deformations \`a la \cite{fmn}, thus giving rise to integrable dispersive (3+1)-dimensional systems; in this connection cf.\ also e.g.\ \cite{lmp} and references therein on quantization of the contact bracket.

Next, as the nonlinear contact Lax pairs
do not involve variable spectral parameter, they are likely to lend themselves to
discretization considerably easier than linear contact Lax pairs.

In turn, this should make the construction of integrability-preserving discretizations (cf.\ e.g.\ \cite{b-s} and references therein on the latter) of integrable systems (\ref{sys}) associated with (\ref{nLax4D})
far easier than that of other integrable (3+1)-dimensional dispersionless systems known to date.

On the other hand,
the potential availability of plethora of exact solutions for integrable systems associated with (\ref{nLax4D})
via the inverse scattering transform for (\ref{clp}) or other methods, as discussed in more detail in Section~\ref{lin-lax}, should make such systems {\em inter alia}
useful benchmarks for the known and new numerical methods of solving general
(3+1)-dimensional dispersionless systems. 
\looseness=-1


Notice that in a degenerate special case when $f$ and $g$ are affine in $p$,
system (\ref{nLax4D}) becomes a linear system for $\psi$ of the form
\begin{equation}\label{linfg}
\psi_y=f_0\psi_z+f_1\psi_x,\quad \psi_t=g_0\psi_z+g_1\psi_z,
\end{equation}
where $f_0,f_1,g_0,g_1$ are functions of $\bu$.

If $f_0,f_1,g_0,g_1$ further involve an additional nonremovable parameter, say, $\lambda$,
then (\ref{linfg}) becomes a linear isospectral Lax pair with the spectral parameter $\lambda$ while (\ref{clp})
becomes irrelevant.

There exist, see e.g.\ \cite{d-book,fkk,km,ms,zs} and references therein,
integrable (3+1)-dimensional systems with the Lax pairs of this form, i.e.,
of the form
(\ref{linfg}) with $f_0,f_1,g_0,g_1$ depending on $\bu$ and $\lambda$; perhaps the best known examples of this kind are provided by the heavenly equations,  see e.g.\ \cite{d-book}. These systems fit into our approach
as (degenerate) special cases, but everywhere else in the present paper we concentrate on the case when the Lax functions $f(p,\bu)$ and
$g(p,\bu)$ are inherently nonlinear in $p$. 

In closing 
let us point out a more symmetric parametric form of the nonlinear contact Lax pair (cf.\ e.g.\ \cite{os}
in (2+1) dimensions),
\begin{equation}\label{laxpa}
\psi_x=\psi_z \Phi(\bu,\zeta),\quad
\psi_y=\psi_z\Theta(\bu,\zeta),\quad \psi_t=\psi_z\Omega(\bu,\zeta),
\end{equation}
where $\zeta=\zeta(x,y,z,t)$. 

Expressing $\zeta$ from the first equation of (\ref{laxpa}) and
substituting the result into the remaining two equations gets us
back to (\ref{nLax4D}). On the other hand, (\ref{nLax4D}) can be written as a special case of (\ref{laxpa}):
\[
\psi_x=\psi_z \zeta,\quad
\psi_y=\psi_z F(\bu,\zeta),\quad \psi_t=\psi_z G(\bu,\zeta).
\]


\section{Explicit examples of new (3+1)-dimensional integrable systems}\label{exa}
\subsection*{\em Example 1}
Let $N=4$, $\boldsymbol{u}=(u,v,w,r)^{\mathrm{T}}$,
$f=v p^{-1}+u$, $g=w p+r p^{2}$, 
so the associated linear contact Lax pair (\ref{clp}) reads
\begin{equation}\label{ralinlax}
\begin{array}{l}
\chi_y=\displaystyle\!\!\left(p
u_z-u_x+v_z-\frac{v_x}{p}\right)\!\chi_p-\frac{v}{p^2}\chi_x+\!\left(u+\frac{2
v}{p}\right)\!\chi_z,
\\[5mm]
\chi_t=\displaystyle\!\!p\left(r_z
p^2+(w_z-r_x)p-w_x\right)\!\chi_p+\left(2 p r+w\right)\!\chi_x-r
p^2\chi_z.
\end{array}
\end{equation}

The nonlinear contact Lax pair (\ref{nLax4D}) in this case has the form
\begin{equation}\label{ralax}
\psi_y=v \psi_z^2/\psi_x+u \psi_z,\quad \psi_t=w\psi_x+r\psi_x^2/\psi_z. 
\end{equation}

The compatibility condition (\ref{zcr}) yields, upon equating to
zero the coefficients at all powers of $p$,
a non-evolutionary system
\begin{equation}\label{sysuvwr}
\begin{array}{rcl}
u_t &=& 2 r v_x-2 v w_z+v r_x+w u_x,\\[3mm]
v_t &=& v w_x+w v_x,\\[3mm]
w_y &=& 2 v r_z-2 r u_x+u w_z+r v_z,\\[3mm]
r_y &=& u r_z+r u_z.
\end{array}
\end{equation}

Note that systems (\ref{ralinlax}), (\ref{ralax}) and
(\ref{sysuvwr}) are invariant under the simultaneous swap of the
following pairs of variables: $y\leftrightarrow t$,
$x\leftrightarrow z$, $u\leftrightarrow r$, $v\leftrightarrow w$.


The second and fourth equations of (\ref{sysuvwr}) have the form of conservation laws, so
we can
introduce the potentials $a$ and $b$
such that $v=a_x$, $v w =a_t$, $r=b_z$, $u r=b_y$, whence
\[
u=b_y/b_z,\quad v=a_x,\quad w=a_t/a_x,\quad r=b_z. 
\]
This 
turns (\ref{sysuvwr}) into a system of two coupled second-order equations for $a$ and $b$,\looseness=-1
\begin{equation}\label{sysab}
\begin{array}{rcl}
a_{yt} &=& 2 a_x^2 b_{zz}-2 a_x
b_{xy}+\displaystyle\frac{a_t}{a_x}a_{xy}+\frac{b_y}{b_z}a_{zt}\\[4mm]&&+\displaystyle\frac{2
a_x b_y}{b_z}b_{xz}
+\frac{(a_x^2 b_z^2-a_t b_y)}{a_x b_z}a_{xz},\\[5mm]
b_{yt} &=& 2 b_z^2 a_{xx}-2 b_z
a_{zt}+\displaystyle\frac{b_y}{b_z}b_{zt}+\frac{a_t}{a_x}b_{xy}\\[4mm]&&+\displaystyle\frac{2
a_t b_z}{a_x}a_{xz} +\frac{(a_x^2 b_z^2-a_t b_y)}{a_x b_z}b_{xz},
\end{array}
\end{equation}
which also inherits the above discrete symmetry: it is invariant
under the simultaneous swap $y\leftrightarrow t$, $x\leftrightarrow
z$, $a\leftrightarrow b$.

Plugging a formal expansion
$\chi=\sum\limits_{\mu=0}^\infty\omega_\mu p^{\mu}$ into (\ref{ralinlax})
yields
\[
\begin{array}{rcl}
\hspace*{-9mm}(\omega_\mu)_x\!\!&=&\!\!v^{-1}\left((\mu-2)
u_z\omega_{\mu-2}-(\mu-1)(u_x-v_z)\omega_{\mu-1}-\mu v_x \omega_\mu \right.\\[2mm]
&&\left.
-(\omega_{\mu-2})_y+ u(\omega_{\mu-2})_z+2 v(\omega_{\mu-1})_z\right),\\[3mm]
\hspace*{-9mm}(\omega_\mu)_t\!\!&=&\!\!(\mu-2)r_z \omega_{\mu-2}-r
(\omega_{\mu-2})_z+ (\mu-1)
(w_z-r_x) \omega_{\mu-1}\\[2mm]
&&-j w_x\omega_{\mu}+2 r (\omega_{\mu-1})_x+w(\omega_{\mu})_x,
\end{array}
\]
where $\mu=0,1,2,\dots$. 
Here and below we tacitly assume that 
$\omega_\mu=0$ for
$\mu<0$.\looseness=-1

For $\mu=0$ we have $(\omega_0)_x=(\omega_0)_t=0$, so
$\omega_0=\omega_0(y,z)$, i.e., $\omega_0$ is not a nontrivial
nonlocal variable but 
just an arbitrary
(sufficiently smooth) function of $y$ and $z$.

Substituting this into the compatibility conditions $((\omega_\mu)_x)_t=((\omega_\mu)_t)_x$ for $\mu=1,2,\dots$, yields,
in a slight modification of the setup discussed
after (\ref{col}),
an infinite hierarchy of nonlocal conservation laws for (\ref{sysuvwr}), where $\mu=1,2,\dots$:
\[
\begin{array}{l}
\left(v^{-1}\left((\mu-2) u_z\omega_{\mu-2}-(\mu-1)(u_x-v_z)\omega_{\mu-1}\right.\right.\\[2mm]\left.\left.-\mu v_x \omega_\mu-(\omega_{\mu-2})_y+ u(\omega_{\mu-2})_z+2 v(\omega_{\mu-1})_z\right)\right)_t\\[2mm]
=\left((\mu-2)r_z \omega_{\mu-2}+ (\mu-1) (w_z-r_x)
\omega_{\mu-1}\right.\\[2mm]\left.-j w_x\omega_{\mu}+2 r (\omega_{\mu-1})_x+w(\omega_{\mu})_x -r
(\omega_{\mu-2})_z\right)_x\!.
\end{array}
\]
We conjecture that all these nonlocal conservation laws are nontrivial.


So far we were unable to find a Hamiltonian structure for (\ref{sysuvwr}), or, more precisely, a generalized Hamiltonian structure in the sense of \cite{kv} and references therein, as (\ref{sysuvwr}) is not in evolutionary form,
thus leaving open the matters of finding an interpretation for the above conservation laws as well as of the study of
the (generalized) Poisson brackets of the associated functionals.

\subsection*{\em Example 2}
Consider the special case of (\ref{ratfg})
when $m=1$ and $n=2$.
To simplify notation, 
put $v_0=u$, $v_1=w$, $w_0=v$, $w_1=r$, so
\[
f=p^2+w p+u,\quad g=p^3+2 w p^2+r p+v.
\]
Thus, the nonlinear contact Lax pair reads
\[
\begin{array}{rcl}
\displaystyle\psi_y&=&\psi_z \left(\left(\displaystyle\frac{\psi_x}{\psi_z}\right)^2 +w\displaystyle \frac{\psi_x}{\psi_z}+u \right),\\[7mm] \displaystyle \psi_t&=&\psi_z \left(\left(\displaystyle\frac{\psi_x}{\psi_z}\right)^3+2 w \left(\displaystyle\frac{\psi_x}{\psi_z}\right)^2+r\displaystyle \frac{\psi_x}{\psi_z}+v\right),
\end{array}
\]
while the linear contact Lax pair takes the form
\[
\begin{array}{rcl}
\chi_y&=&\left(2p+w\right)\chi_x+\left(u-p^2\right)\chi_z
+\left(w_z p^2+(u_z-w_x) p-u_x\right)\chi_p,\\[5mm]
\chi_t&=&\left(3 p^2+4 w p+r\right)\chi_x+\left(-2 p^3-2 w p^2+v\right)\chi_z\\[5mm] &&+\left( 2 w_z p^3+(-2 w_x+r_z) p^2+(-r_x+v_z) p-v_x)\right)\chi_p.
\end{array}
\]

The compatibility condition for either of the above two Lax pairs has the form of a four-component system
\begin{equation}\label{4ddkp}
\begin{array}{l}
u_t-v_y-v u_z-r u_x+u v_z+w v_x=0,\\[3mm]
2 u_z-r_z+w_x+2 w w_z=0,\\[3mm]
2 r_x-3 u_x-2 w_y-v_z+2 w u_z-2 w w_x+2 u w_z=0,\\[3mm]
w_t-r_y+2 v_x-4 w u_x+w r_x-r w_x-v w_z+u r_z=0,
\end{array}
\end{equation}
which, in perfect agreement with the general results from Section~\ref{mr},
can be solved with respect to $u_z,v_z,w_z,r_z$, so we have an evolution system in disguise.
\looseness=-1

The second equation in (\ref{4ddkp}) has the form of a conservation
law and allows us, by analogy with Example~1, to introduce a
potential $a$ for this conservation law such that $a_z=w$ and
$a_x=r-w^2-2u$, whence
\[
\begin{array}{l}
w=a_z, 
\\[3mm]
r=a_z^2+2 u+a_x.
\end{array}
\]
Plugging these expressions into the third equation of (\ref{4ddkp})
we can, in turn, rewrite this equation as a conservation law and
introduce a potential $b$ for this conservation law so that
\[
\begin{array}{l}
u=-a_z^2-2 a_x-b_z,
\\[3mm]
v=-2 a_y-2 a_z^3-4 a_z a_x-2 a_z b_z-b_x.
\end{array}
\]

Substituting the resulting expressions into the remaining
equations of (\ref{4ddkp}) yields the following second-order system
for $a$ and $b$:
\begin{equation}\label{4dpdkp}
\begin{array}{rcl}
a_{zt}&=&3 a_z a_{xx} + a_{xy} + (2 a_z^2 - a_x - b_z) a_{xz} - 2 a_z a_{yz}\\[1mm] && - (4 a_z^3 + 8 a_x a_z + 4 a_z b_z + 2 a_y + b_x) a_{zz}\\[1mm] && + 2 b_{xx}+ 2 a_z b_{xz} - 2 b_{yz} - (2 a_z^2 + 4 a_x + 2 b_z) b_{zz},\\[2mm]
b_{zt}&=&-(12 a_z^2 + 6 a_x + 4 b_z) a_{xx} + 2 a_{yy} + (12 a_z^2 + 8 a_x + 4 b_z) a_{yz}\\[1mm] &&- (12 a_z^3 +8 a_x a_z + 4 a_z b_z + 4 a_y + 2 b_x) a_{xz}\\[1mm] && + 2 (5 a_z^2 + 2 a_x + b_z)(a_z^2 + 2 a_x + b_z) a_{zz} - 2 a_{xt} - 5 a_z b_{xx}\\[1mm] && + b_{xy} -(6 a_z^2 + a_x +b_z) b_{xz} + 6 a_z b_{yz}\\[1mm] && + (4 a_z^3 + 8 a_x a_z + 4 a_z b_z - 2 a_y - b_x) b_{zz}.
\end{array}
\end{equation}

If $a$ and $b$ are independent of $z$ then this system boils down to
\[
a_{xy}+2 b_{xx}=0,\quad b_{xy}-6 a_x a_{xx} + 2 a_{yy} - 2 a_{xt}=0,
\]
and upon imposing a reduction
$a=s_x, b=-s_y/2$
we end up with a single
equation
\[
12 s_{xx}s_{xxx}-3 s_{xyy}+4 s_{xxt}=0,\vspace{-2mm}
\]
or, in terms of $a$,
\[
12 a_x a_{xx}-3 a_{yy}+4 a_{xt}=0. 
\]

Up to a rescaling of $x, y$, $t$, and $a$, the above equation is
nothing but the potential dispersionless Kadomtsev--Petviashvili
(dKP) equation, cf.\ e.g.\ \cite{d-book,kg2,imk,z}, also known as
the (2+1)-dimensional potential Khokhlov--Za\-bo\-lot\-skaya
\cite{zakh} equation or as the (2+1)-dimensional
Lin--Reissner--Tsien \cite{lrt} equation.\looseness=-1

Thus, (\ref{4dpdkp}), or, equivalently, (\ref{4ddkp}), provides a
novel (3+1)-di\-men\-si\-on\-al {\em integrable}
generalization of the potential dKP equation, unlike, say, the
nonintegrable (3+1)-di\-men\-si\-on\-al potential
Khokhlov--Za\-bo\-lot\-skaya equation, cf.\ e.g.\ \cite{ros},
which,
again up to a suitable rescaling of $x,y,z,t,a$, reads
\[ \left(4 a_t+6 a_x^2\right)_x- 3 (a_{yy}+a_{zz})=0. \]

\subsection*{Acknowledgments}
The author is pleased to thank M. B\l aszak, 
I.S. Krasil'shchik, M. Kunzinger, S. Leble, B. McKay, O.I. Morozov, P.J. Olver, R.O. Popovych, V. Rubtsov, I.A.B. Strachan, L. Vitagliano, and R. Vitolo for stimulating discussions and helpful comments, and to the anonymous referees for useful suggestions.



\end{document}